\documentclass{article}
\usepackage{microtype}
\usepackage{verbatim}
\usepackage{graphicx}
\usepackage{rotating}
\usepackage{listings}

\usepackage{cmll}

\usepackage{url}
\usepackage{hyperref}

\usepackage{pslatex}

\usepackage{datetime}
\newdateformat{versiondate}{%
\THEMONTH\THEDAY}

\usepackage{pgf,tikz,environ}
\usetikzlibrary{arrows}
\usetikzlibrary{cd}

\usepackage{xpatch} 

\usepackage{amsthm}
\usepackage[fleqn]{amsmath} % flush equations left
\usepackage{amsfonts}
\usepackage{amssymb}

\newtheoremstyle{zoltanstyle}
  {1em} % Space above
  {\topsep} % Space below, usually \topsep
  {} % Body font
  {} % Indent amount
  {\bfseries} % Theorem head font \bfseries
  {.} % Punctuation after theorem head
  {.5em} % Space after theorem head
  {} % Theorem head spec (can be left empty, meaning `normal')

\theoremstyle{zoltanstyle}
\swapnumbers
\xpatchcmd\swappedhead{~}{.~}{}{}
\newtheorem{body}{}
\numberwithin{body}{section}

\newtheorem{corollary}[body]{Corollary}
\newtheorem{definition}[body]{Definition}
\newtheorem{example}[body]{Example}

\newtheorem{lemma}[body]{Lemma}

\newtheorem{proposition}[body]{Proposition}

\newtheorem{theorem}[body]{Theorem}

\expandafter\let\expandafter\oldproof\csname\string\proof\endcsname
\let\oldendproof\endproof

\renewenvironment{proof}[1][\proofname]{%
  \oldproof[\normalfont \bfseries #1.]%
}{\oldendproof}

\usepackage{enumitem}
\usepackage{mathptmx}
\usepackage{xspace}

\usepackage[T1]{fontenc}
\usepackage[utf8]{inputenc}

\newcommand{\pita}{\mathrel{\reflectbox{\rotatebox[origin=c]{180}{$\mathbb{A}$}}\!\!}}  % PITts forAll
\newcommand{\pite}{\mathrel{\reflectbox{$\mathbb{E}$}\!\!}}  % PITts Exists
\newcommand{\vdashl}{\vdash_{\mathcal{L}}}
\newcommand{\vdashat}{\vdash_{at}}

\title{Proof-theoretic methods in\\ quantifier-free definability}
\author{Zoltan A. Kocsis}
\date{21 January 2025}

\begin{document}

\maketitle

\begin{abstract}
We introduce a proof-theoretic approach to showing nondefinability of second-order intuitionistic connectives by quantifier-free schemata. We apply the method to prove that Taranovsky's "realizability disjunction" connective does not admit a quantifier-free definition, and use it to obtain new results and more nuanced information about the nondefinability of Kreisel’s and Po\l acik’s unary connectives. The finitary and combinatorial nature of our method makes it resilient to changes in metatheory, and suitable even for settings with axioms that are explicitly incompatible with classical logic. Furthermore, the problem-specific subproofs arising from this approach can be readily transcribed into univalent type theory and verified using the Agda proof assistant.
\end{abstract}

\section{Introduction} \label{sec:introduction}

\begin{body}
We present a new, purely proof-theoretic method for showing undefinability of second-order intuitionistic connectives by quantifier-free schemata. This new approach leverages the Pitts quantifier theorem, a tool that has existed for three decades but hitherto had limited application in this setting, chiefly due to the computational difficulty involved in manually calculating the required formulae. However, the year 2023 marked a pivotal development: F{\'{e}}r{\'{e}}e and Van Gool~\cite{feree-pitts} created a formally verified implementation of the algorithm inherent in the proof of the Pitts quantifier theorem. They not only give a formalization of Pitts' result in the Coq proof assistant, but in addition provide a reasonably efficient, correct-by-construction OCaml program extracted from this Coq formalization, which can compute Pitts interpolants for arbitrary propositional formulae. As we shall see, this paves the way to rapid, streamlined, robust proofs of non-definability for second-order connectives.
\end{body}

\begin{body}
Given a connective $\varphi$ definable in second-order intuitionistic propositional logic, one can summarize the method presented in the article using the following three steps:
\begin{enumerate}
\item We compute the Pitts interpolants $\pita X. \varphi$ and $\pite X. \varphi$.
\item Conditional on the definability of the connective by a quantifier-free schema, we use the information obtained from the Pitts interpolants to deduce the existence of finitely many auxiliary formulae $t_1,t_2,\dots t_n$ that witness the definability in a certain sense.
\item We show, using the consequence relation of intuitionistic logic, that the existence of the formulae $t_1,\dots,t_n$ is sufficient to deduce a formula that is not intuitionistically valid (usually the law of excluded middle), thereby showing that the terms $t_1,\dots,t_n$ cannot exist in intuitionistic logic after all.
\end{enumerate}
Thanks to the formally verified implementation of the Pitts quantifier theorem, the first step is entirely automated. The second step of the process is largely mechanical as well: given the Pitts interpolants, the formulae in question can be obtained by analysis of cut-free proofs of the definability of the connective. The third step does require ingenuity: one has to find the appropriate intuitionistically invalid formula, then conjure up a proof of it by reasoning purely syntactically, in terms of derivable sequents. The difficulty of this step is comparable to choosing an appropriate topological space in a non-definability argument based on topological semantics.
\end{body}

\begin{body}
The finitary and combinatorial nature of the method presented here makes it more resilient against metamathematical alterations, such as passing to a constructive metatheory (or a metatheory which is not compatible with classical logic), than the widely applied strategies based on topological semantics. Additionally, the subproofs which emerge during the method's third step have a level of generality that allows for direct reinterpretation within univalent type theories, such as Homotopy Type Theory. This not only gives us new type-theoretic results at no additional conceptual cost, but also allows us to check our proofs by reformulating them in the language of an interactive theorem prover based on type theory. Indeed, all results of this sort that we present in the article have been formalized and checked using the Agda proof assistant.
\end{body}

\begin{body} \textbf{Background.}
Questions about quantifier-free definability of second-order intuitionistic connectives were first raised by \textit{Kreisel}~\cite{kreisel-monadic}. \textit{Troelstra}~\cite{troelstra-kreiselstar} used topological semantics to characterize definability of certain families of monadic second-order connectives. Later, \textit{Po\l acik}~\cite{polacik-pitts} introduced an undefinable monadic connective whose undefinability cannot be shown via topological semantics in any dense-in-itself metric space. Our method not only gives syntactic proofs of both undefinability results, but strengthens them by extracting information about which super-intuitionistic logics \textit{can} provide quantifier-free definitions of similar connectives. Beyond these classic works, there are several open questions of interest concerning quantifier-free definability. One such question concerns a  disjunction-like connective proposed by Taranovsky and motivated by realizability considerations. Bauer~\cite{bauer-equivalence} showed that this connective is definable using propositional quantifiers, and asked whether one can find an equivalent quantifier-free definition. This article gives a negative answer.
\end{body}

\begin{body} \textbf{Outline.}
We introduce the necessary conventions and notations in the remainder of Section~\ref{sec:introduction}, and explain how the formalization of the Pitts quantifier theorem enables our work. The general results required to apply the new method, most importantly Theorem~\ref{thm:general:auxiliaryequiv}, are established in Section~\ref{sec:general}. The next three sections apply our method to interesting families of binary and unary connectives. Section~\ref{sec:tara} treats Taranovsky's realizability-inspired disjunction connective: the main result is Theorem~\ref{thm:tara:mainresult}, which shows that this connective is not definable by a quantifier-free schema in intuitionistic logic. Section~\ref{sec:kreisel} focuses on Kreisel's star connective, and a family of unary connectives related to it. Theorem~\ref{thm:kreisel:troelstra} re-proves a classical result of Troelstra stating that none of the connectives in the family are definable by quantifier-free schemata. Analyzing the new proof, we obtain new results which give finer-grained information about which logics can and cannot define these connectives (Corollaries~\ref{cor:kreisel:extensions} and \ref{cor:kreisel:disjunction}). Section~\ref{sec:polacik} treats a connective introduced by Po\l acik. We re-prove the nondefinability of Po\l acik's connective in intuitionistic logic (Theorem~\ref{thm:polacik:mainresult}) using proof-theoretic methods, obtain new nondefinability results for closely related connectives (Corollary~\ref{cor:polacik:wlem}) and give a condition on logics which do admit a quantifier-free definition (Corollary~\ref{cor:polacik:disjunction}). The background on type-theoretic results and their Agda formalizations is summarized in Section~\ref{sec:typetheory}, while the individual results are stated in the sections of the respective connectives.
\end{body}

\subsection{Second-order propositional logic}

\begin{definition}
The \textit{language of intuitionistic second-order propositional logic} consists of
\begin{itemize}
\item \textbf{Propositional variables}: countably many propositional variable symbols, usually denoted by upper-case letters $A,B,C,\dots,P,Q,R,\dots,X, Y, Z$;
\item \textbf{Nullary connectives}: the falsum constant denoted $\bot$;
\item \textbf{Binary connectives}: conjunction denoted $\wedge$, disjunction denoted $\vee$, implication denoted $\rightarrow$; and
\item \textbf{Propositional quantifiers}: existential denoted $\exists X$ and universal denoted $\forall X$ where $X$ stands for some propositional variable symbol.
\end{itemize}
We define \textit{formulae}, the notion of \textit{free and bound occurrences} of variables in formulae, and \textit{simultaneous substitution} of terms $T_1,\dots,T_n$ for distinct variables variable $X_1,\dots,X_n$ (denoted $A[T_1/X_1,\dots,T_n/X_n]$) in a capture-avoiding manner in the obvious way.
\end{definition}

\begin{body}
Although the connectives $\wedge,\vee,\bot$ and even the quantifier $\exists$ admit definitions in terms of $\rightarrow$ and $\forall$ in second-order logic (see~\cite{urzyczyn-lectures}),  our language for the logic incorporates them as primitives instead. This choice stems from the fact that the definitions of these connectives inherently involve universal quantifiers, and we aim for the quantifier-free fragment of our second-order logic to correspond exactly to ordinary (non-second-order) intuitionistic propositional logic.  We do introduce two abbreviations, however: $\neg X$ standing for $X \rightarrow \bot$,  and $X \leftrightarrow Y$ standing for $(X\rightarrow Y) \wedge (Y \rightarrow X)$.
\end{body}

\begin{definition}\label{def:general:sequentrules}
As usual, a \textit{sequent} is an expression $\Gamma \vdash A$ where $A$ is a formula and $\Gamma$ denotes a finite (possibly empty) sequence of formulae, considered up to order. A sequent $\Gamma$ with empty right-hand side denotes $\Gamma \vdash \bot$. We follow the usual presentation in the style of Gentzen and Girard, and take the following \textit{inference rules for second-order propositional logic}: \\
\begin{tabular}{lllll}
 \textbf{Identity} & & & \\
 & $\frac{~}{\varphi \vdash \varphi}\text{ax}$ & $\frac{\Gamma \vdash \varphi ~~~~ \Delta, \varphi \vdash \psi}{\Gamma, \Delta \vdash \psi}\text{cut}$ \\ \\
 \textbf{Structure} & & & \\
 & $\frac{\Gamma \vdash \psi}{\Gamma, \varphi \vdash \psi}wL$ & $\frac{\Gamma, \varphi, \varphi \vdash \psi}{\Gamma, \varphi \vdash \psi}cL$ & $\frac{\Gamma \vdash}{\Gamma \vdash \varphi}wR$ \\ \\
 \textbf{Connectives} & & & \\
 & % Disj
 $\frac{\Gamma, \varphi \vdash \theta ~~~~ \Gamma, \psi \vdash \theta}{\Gamma, \varphi \vee \psi \vdash \theta}\vee L$ &
 $\frac{\Gamma \vdash \varphi}{\Gamma \vdash \varphi \vee \psi}\vee R_1$ &
 $\frac{\Gamma \vdash \psi}{\Gamma \vdash \varphi \vee \psi}\vee R_2$
 \\ \\
 & % Conj
 $\frac{\Gamma \vdash \varphi ~~~~ \Gamma \vdash \psi}{\Gamma \vdash \varphi \wedge \psi}\wedge R$ &
 $\frac{\Gamma, \varphi \vdash \theta}{\Gamma, \varphi \wedge \psi \vdash \theta}\wedge L_1$ &
 $\frac{\Gamma, \psi \vdash \theta}{\Gamma, \varphi \wedge \psi \vdash \theta}\wedge L_2$
 \\ \\
 & % Impl
 $\frac{\Gamma \vdash \varphi ~~~~ \Delta, \psi \vdash \theta}{\Gamma, \Delta, \varphi \rightarrow \psi \vdash \theta}\rightarrow L$ &
 $\frac{\Gamma, \varphi \vdash \psi}{\Gamma \vdash \varphi \rightarrow \psi}\rightarrow R$ &
 $\frac{~}{\bot \vdash \varphi}\bot L$
 \\ \\
 \textbf{Quantifiers} & & & \\
 & % Univ
 $\frac{\Gamma, \varphi[T/X] \vdash \psi}{\Gamma, \forall X. \varphi \vdash \psi}\forall L$ &
 $\frac{\Gamma \vdash \varphi}{\Gamma \vdash \forall X. \varphi}\forall R$
 \\ \\
 & % Exis
 $\frac{\Gamma \vdash \varphi[T/X]}{\Gamma \vdash \exists x. \varphi}\exists R$ &
 $\frac{\Gamma, \varphi \vdash \psi}{\Gamma, \exists X. \varphi \vdash \psi}\exists L$
 \\ \\
\end{tabular}
~ \\
where $\Gamma, \Delta$ stand for arbitrary sequences of formulae, $\varphi,\psi,\theta$ stand for arbitrary formulae, and $X$ stands for a propositional variable symbol. As usual, the rules for the quantifiers require some restrictions and clarifications: in the rules $\forall R$ and $\exists L$, the term $T$ should not use bound variables of $\varphi$, while in the rules $\forall L$ and $\exists R$, the variable $X$ must not occur free in the context $\Gamma, \psi$.
\end{definition}

\begin{body}
The rules of Definition~\ref{def:general:sequentrules}, including the cut rule, will feature heavily in the proofs of later sections. We will also make use of the redundancy of the cut rule, and the convergence of the cut-elimination algorithm for second-order propositional logic. For the details of this procedure, and a proof-theoretic argument establishing its convergence, we refer the reader to Girard's writings on the subject~\cite{girard-candred}. Chapter 14 of \textit{Proofs and Types}~\cite{girard-prot} summarizes the key ideas of the proof in terms of System F, the Curry-Howard counterpart to second-order propositional logic.
\end{body}

\begin{body}
In what follows, we will consider derivability both in ordinary (i.e. non-second-order) intuitionistic propositional logic, and in several systems of propositional logic extending intuitionistic logic. Thanks to cut-elimination, and the careful formulation of the language and rules of our system, we need not distinguish between the ordinary and the second-order system: derivability in the quantifier-free fragment of the latter corresponds exactly to the former. Accordingly, we will use the $\vdash$ symbol for both. For derivability in other systems, we will use decorated turnstiles: $\Gamma \vdashl A$ means that the sequent $\Gamma \vdash A$ holds in the logic $\mathcal{L}$. Keep in mind that the observation about the quantifier-free fragment of $\vdash$ need not extend to $\vdashl$.
\end{body}

\begin{definition}\label{def:general:superintuitionistic}
In our text, a \textit{super-intuitionistic (second-order, propositional) logic} $\mathcal{L}$ stands for a binary relation $\vdashl$ between finite sequences of formulae of second-order intuitionistic logic on the left, and a single such formula on the right, such that $\vdashl$ satisfies all the inference rules given in Definition~\ref{def:general:sequentrules}. The \textit{quantifier-free fragment} of such a logic is defined the obvious way.
\end{definition}

\begin{body}\label{sec:superint}
One can obtain a super-intuitionistic logic $\mathcal{L}$ in the sense of Definition~\ref{def:general:superintuitionistic} corresponding to any finitely axiomatizable intermediate propositional logic $L$, simply by defining $\Gamma \vdashl A$ as $\Gamma, \Lambda \vdash A$ where $\Lambda$ consists of the second-order universal closures of the axioms of $L$. In the rest of the article, we use the term \textit{intermediate logic} to refer to systems of this specific form. We will have the opportunity to consider derivability in (the quantifier-free fragments of) several such logics:
\begin{itemize}
\item \textbf{LK}: \textit{classical logic}, given by e.g. the axiom \\
$\forall X. X \vee \neg X$ or $\forall X. \neg\neg X \rightarrow X$
\item \textbf{KC}: the logic of the \textit{weak excluded middle}, given by the axiom \\
$\forall X. \neg X \vee \neg\neg X$;
\item \textbf{KP}: the \textit{Kreisel-Putnam logic} given by the axiom \\
$\forall X. \forall Y. \forall Z. (\neg X \rightarrow (Y \vee Z)) \rightarrow ((\neg X \rightarrow Y) \vee (\neg X \rightarrow Z))$,
\item \textbf{SL}: \textit{Scott's logic}, given by the axiom\\ $\forall X. ((\neg\neg X \rightarrow X) \rightarrow (X \vee \neg X)) \rightarrow (\neg\neg X \vee \neg X)$.
\end{itemize}
\end{body}

\begin{body}
The logics named in \ref{sec:superint} are among the seven distinguished logics between intuitionistic and classical logic which enjoy the Craig~interpolation property~\cite{maksimova-interpolation}, while \textbf{KP} and \textbf{SL} are among the continuum-many~\cite{wronski-disjunctionprop} such systems that enjoy the disjunction property. For in-depth discussion of these systems, we refer the interested reader to Chapter 2 of Fiorentini's thesis~\cite{fiorentini-intermediate}.
\end{body}

\begin{lemma}\label{lemma:general:extensionality}
Assume that the free variables of the formula $P$ are disjoint from the bound variables of the formula $x(P)$. Then the sequent $P \rightarrow P', P' \rightarrow P, x(P) \vdash x(P')$ is provable.
\begin{proof}
Induction on the structure of the formula $x(P)$.
\end{proof}
\end{lemma}

\subsection{Pitts quantifier theorem}

\begin{body}
In the classical version of second-order propositional logic, the second-order quantifiers play no essential role: one can define $\forall X. \varphi$ by means of the formula $\varphi[\bot/X] \wedge \varphi[\neg\bot/X]$, and similarly $\exists X. \varphi$ as $\varphi[\bot/X] \vee \varphi[\neg\bot/X]$. It has been known since the late 1970s that one cannot give such a ``quantifier elimination'' procedure for the intuitionistic variant; and indeed, many of the results in the present article also double as self-contained proofs of this fact.
\end{body}

\begin{body}
As Gabbay~\cite{gabbay-newconnectives} first observed, one can regard the fact that intuitionistic logic does not eliminate its second-order quantifiers as a form of ``incompleteness'' of $\vee, \wedge, \rightarrow, \bot$ \textit{qua} basis for intuitionistic connectives. If we pick a formula $\varphi$ with free variables $P,Q$ which contains one or more quantifiers that cannot be eliminated, we get an altogether new binary connective $\varphi(P,Q)$, one that we cannot express intuitionistically as a compound formula made up of $P,Q,\vee, \wedge, \rightarrow, \bot$, but one which nonetheless remains \textit{classically} equivalent to such a formula.
\end{body}

\begin{body}
Many, but not all quantified formulae give rise to new intuitionistic connectives. For example, the binary connective $\varphi(P,Q)$ given by the formula $\forall Y. (P \rightarrow Y) \rightarrow (Q \rightarrow Y) \rightarrow Y$ is logically equivalent to the ordinary disjunction connective $P \vee Q$ in intuitionistic propositional logic, while the unary connective $\exists Y. X \leftrightarrow \neg Y$ is equivalent to $\neg\neg X \rightarrow X$. The question of which quantified formulae give rise to genuinely new intuitionistic connectives (i.e. the question of definability of connectives by quantifier-free schemata) is the focus of the present article. This question is widely studied for many specific connectives: we will briefly review the prior work on each as they first come up.
\end{body}

\begin{body}
Although one cannot eliminate second-order quantifiers in intuitionistic second-order propositional logic, a celebrated result of Pitts~\cite{pitts-quantifiers} shows that, in a much more limited sense, one can nonetheless model quantification over propositional variables inside the quantifier-free fragment: we have an effectively computable translation $(-)^p$ from the formulae of the full second-order calculus to its quantifier-free fragment which restricts to the identity over quantifier-free formulae, and which is sound in the sense that $\Gamma \vdash \varphi$ implies $\Gamma^p \vdash \varphi^p$.
\end{body}

\begin{body}
The result known variously as the Pitts Quantifier Theorem or the Uniform Interpolation Theorem forms a key ingredient of the aforementioned translation. Here, we state its existential version as Theorem~\ref{thm:general:pitts}: we will not have direct need for the universal version, so we leave that implicit.
\end{body}

\begin{theorem}[Pitts \cite{pitts-quantifiers}]\label{thm:general:pitts}
Consider a finite sequence of propositional variables $\overline{X}$, and a propositional variable $Y$ outside the sequence $\overline{X}$. Let $\Phi(\overline{X}, Y)$ denote a quantifier-free formula containing only the variables in $\overline{X}, Y$. Then one can find a quantifier-free formula $\pite Y.\Phi(\overline{X}, Y)$ so that the following hold:
\begin{enumerate}
    \item All propositional variables in the formula $\pite Y.\Phi(\overline{X}, Y)$ belong to the sequence $\overline{X}$.
    \item For a quantifier-free formula $\Psi(\overline{X})$, intuitionistic propositional logic proves that
    $(\pite Y. \Phi(\overline{X},Y)) \vdash \Psi(\overline{X})$ precisely if it proves that $\Phi(\overline{X},Y) \vdash \Psi(\overline{X})$.
\end{enumerate}
\end{theorem}

\begin{body}
We call the formulae $\pite Y. \varphi$ and $\pita Y. \varphi$ whose existence is asserted by Theorem~\ref{thm:general:pitts} the \textit{Pitts interpolants} of $\varphi$. The proof of the Pitts quantifier theorem is entirely proof-theoretic: it consists of an algorithm that computes the interpolants, and a massive case analysis over proof-trees in Dyckhoff's LJT calculus~\cite{dyckhoff-ljt} that establishes its correctness.
\end{body}

\begin{body}
One can hardly overstate the theoretical significance of Pitts' result: various applications and follow-up works emerged not just in proof theory, but in algebra~\cite{vangool-algebra,kocsis-degreesat} and even classical model theory~\cite{ghilardi-modelcompletions} as well. Alternative, semantic proofs of the result have also become available since its original publication~\cite{ghilardi-sheaves,visser-bisimulations}. However, in practical situations, where one has a specific second-order formula to investigate, the situation was not nearly as rosy. Pitts' algorithm is intricate to an extreme degree, proceeding via lengthy lookup tables, and requiring the mathematician to keep track of large amounts of partial data along the way. The number of recursive calls is also exponential in the nesting depth of implication connectives of the input formula. These factors make the algorithm nigh-impossible to execute by hand for any substantial input, and even difficult to implement correctly on a computer: to the best of our knowledge, no complete implementation was available for the first 30 years of the algorithm's existence.
\end{body}

\begin{body}
Recently, F{\'{e}}r{\'{e}}e and Van Gool~\cite{feree-pitts} gave a formally verified proof of the Pitts quantifier theorem using the Coq proof assistant. Not restricted to a mere formalization of Pitts' result, their work also provides an unexpectedly efficient, correct-by-construction OCaml program (\texttt{propquant}), derived directly from the Coq formalization, which is able to calculate the Pitts interpolants of any given formula. This opens the way towards many approaches to proof-theoretic problems which would have been infeasible without a reliable way of calculating interpolants. As the following sections will elucidate, it in particular paves the way toward streamlined and robust proofs of non-definability for a large number of second-order connectives.
\end{body}

\subsection{Semantic matters}

\begin{body}
When treating quantifier-free definability, prior works almost universally rely on one of several sound semantics for intuitionistic second-order propositional logic. If one can exhibit a semantic model in which a formula with quantifiers has a different denotation than any quantifier-free formula, that suffices to demonstrate the non-definability of the formula in question. The sound semantics deployed for this endeavor can be broadly categorized into several interconnected families: algebraic, Kripke-style, topological, and realizability/topos interpretations. Topological approaches stand out as the prevalent methodologies among these, serving as a fundamental framework for numerous inquiries and analyses in the field. For overviews of these semantic methods, we refer the reader to \cite{urzyczyn-lectures,polacik-pitts,zdanowski-existential}.
\end{body}

\begin{body}
In this work, we pivot in a different direction, using methods anchored not in semantics, but in structural proof theory. We aim to demonstrate, through the applications presented below, that the recent computational advancements in Pitts quantifiers make such methods competitive with the standard semantic approaches. The structural proofs yield new information about definability and non-definability in various logics, work uniformly in a wide variety of metatheories, and their key parts are easy to transcribe inside the internal logic of proof assistants.
\end{body}

\begin{body}
The resilience to variations in the metatheory is a particularly pertinent facet: for philosophical and pragmatic reasons, practitioners will naturally prefer to work within an intuitionistic foundational system while investigating definability in second-order intuitionistic logic. The formulation of general topology within a sufficiently constructive setting comes with its own difficulties, but even specific scenarios pose interesting metatheoretic challenges. Troelstra was the first to highlight the acute sensitivity of topological semantics for second-order logic to the chosen metatheory. A noteworthy instance involves Kreisel's unary star connective $\ast(P)$ (treated in detail in Section~\ref{sec:kreisel}) which coincides with $\neg\neg P$ in topological semantics applied over $\mathbb{R}$ as long as one operates within a classical metatheory. However, inside a constructive metatheory, one has no hope of proving this: it would contradict Church's thesis!\footnote{Certain varieties of constructive mathematics, such as Russian Constructivism, admit Church's thesis as a theorem, and can turn this observation into a semantic proof that $\ast(P)$ has no quantifier-free definition! One cannot translate this proof into, say, CZF, which does not admit Church's thesis as a theorem: there, one needs a different semantic argument using an altogether different space.} This is in strict contrast with the quantifier-free case: by a famous result by McKinsey and Tarski~(see \cite{kremer-strongcomplete}), dense-in-itself metric spaces like $\mathbb{R}$ provide not just sound, but complete semantics in the quantifier-free case. Similar phenomena affect certain nondefinability arguments involving realizability and topoi in constructive metatheories. Without going into details, we remind the reader of the striking fact that one may obtain non-preorder complete small categories when working internally to a topos (see Hyland's~\cite{hyland-smallcomplete} article for an example); working in classical sets, no such categories exist~\cite{shulman-settheory}.
\end{body}

\begin{body}
The dependence on metatheory makes semantic proofs of quantifier-free nondefinability of connectives quite inconvenient: before asserting nondefinability in one's chosen foundational system, one has to check in detail whether the proof in the literature actually applies in the current setting, for example whether the space used by the proof can be constructed at all (similar phenomena arise frequently in the first-order setting as well; e.g. the results of Lubarsky~\cite{lubarsky-topology} show that distinguishing LLPO from LPO using topological models essentially requires a non-principal ultrafilter). In contrast, proof-theoretic methods are combinatorial in nature, and work uniformly in a very wide variety of foundational settings. Essentially, any foundational theory with the power to establish cut-elimination for second-order propositional logic is capable of carrying out these structural arguments in an identical manner. While proving cut-elimination for second-order propositional logic is admittedly a potent property in terms of consistency strength (e.g. Heyting arithmetic cannot prove said result), it is nonetheless fully constructive~\cite{altenkirch-systemf}.
\end{body}

\section{Regular connectives}\label{sec:general}

\begin{body}
In this section we introduce the notion of regular connective (Definition~\ref{def:general:regularconn}), our primary object of study for the rest of the article. We outline the connection between Pitts quantifiers and definability of these connectives in intuitionistic propositional logic. This paves the way for the introduction of auxiliary formulae (Definition~\ref{def:general:auxiliary}), whose analysis plays a key role in deducing definability properties of the specific connectives treated below.
\end{body}

\begin{definition}\label{def:general:regularconn}
Consider a sequence of $n$ distinct propositional variables $X_1,\dots,X_n$, and a propositional variable $Y$ outside this sequence. Let $\Phi(X_1,\dots,X_n, Y)$ denote a quantifier-free formula containing only the variables in $X_1,\dots,X_n,Y$. We introduce the notation $C_{\Phi}(X_1,\dots,X_n)$ for the formula $\exists Y. \Phi(X_1,\dots,X_n, Y)$, and call $C_{\Phi}(X_1,\dots,X_n)$ the \textit{regular connective of arity $n$ defined by $\Phi$}.
\end{definition}

\begin{body}
Generally for a sequence of formulae $\varphi_1,\dots,\varphi_n$, we regard $C_{\Phi}(\varphi_1,\dots \varphi_n)$ as an abbreviation for the formula $$(\exists Y. \Phi(X_1,\dots,X_n, Y))[\varphi_1/X_1,\dots,\varphi_n/X_n].$$ When possible, we leave the arity of regular connectives implicit, and write $\Phi(\overline{X})$ instead of $\Phi(X_1,\dots, X_n)$ to indicate a sequence of variables of the appropriate length.
\end{body}

\begin{body}
In parallel to Definition~\ref{def:general:regularconn}, one could introduce $\forall$-regular connectives that use a universal quantifier in place of the existential one. While many of the observations and results possess analogous counterparts in the universal context, they are not always as straightforward to establish. Since contraction is forbidden on the right, where $\exists$~substitution happens, but not on the left, where $\forall$~substitution happens, most results which rely on analyzing a single formula in the existential setting require analyzing a whole finite sequence of them in the universal setting. While these difficulties are not insurmountable, developing the theory for $\forall$-regular connectives is beyond the scope of this article.
\end{body}

\begin{definition}\label{def:general:definability}
Consider a regular connective $C_\Phi(\overline{X})$ and a quantifier-free formula $\Psi(\overline{X})$ whose variables belong to the same sequence $\overline{X}$. We say that $\Psi$ \textit{defines} $C_\Phi$ \textit{in intuitionistic propositional logic} if
\begin{enumerate}
    \item $C_\Phi(\overline{X}) \vdash \Psi(\overline{X})$ and
    \item $\Psi(\overline{X}) \vdash C_\Phi(\overline{X})$
\end{enumerate}
both hold. Naturally, we call the connective $C_\Phi(\overline{X})$ \textit{definable in intuitionistic propositional logic}, or \textit{definable by a quantifier-free schema} if we can find some quantifier-free formula that defines it.
\end{definition}

\begin{proposition}\label{prop:general:erightrule}
Assume that the propositional variable $X$ does not occur in the contexts $\Gamma, \Delta$. The second-order propositional calculus derives $\Gamma, \exists X. \Phi(X) \vdash \Delta$ precisely if it derives $\Gamma, \Phi(X) \vdash \Delta$.
\begin{proof}
A straightforward commutation argument, using the fact that the second-order propositional calculus has cut-elimination.
\end{proof}
\end{proposition}

\begin{proposition}\label{prop:general:definability}
If a regular connective $C_\Phi(\overline{X})$ is definable in intuitionistic propositional logic, then the formula $\pite Y. \Phi(\overline{X}, Y)$ defines $C_\Phi$.
\begin{proof}
Assume that we can find some quantifier-free formula $\Psi(\overline{X})$ so that the second-order calculus proves both
\begin{enumerate}
    \item $C_\Phi(\overline{X}) \vdash \Psi(\overline{X})$, and
    \item $\Psi(\overline{X}) \vdash C_\Phi(\overline{X})$
\end{enumerate}
It follows from Theorem~\ref{thm:general:pitts} that $C_\Phi(\overline{X}) \vdash \pite Y. \Phi(\overline{X}, Y)$. Cut against $\Psi(\overline{X}) \vdash C_\Phi(\overline{X})$ to obtain $\Psi(\overline{X}) \vdash \pite Y. \Phi(\overline{X}, Y)$. From Proposition~\ref{prop:general:erightrule}, we know that $C_\Phi(\overline{X}) \vdash \Psi(\overline{X})$ implies $\Phi(\overline{X},Y) \vdash \Psi(\overline{X})$. Theorem~\ref{thm:general:pitts} shows that this is equivalent to $\pite Y. \Phi(\overline{X},Y) \vdash \Psi(\overline{X})$. Applying the cut rule, we get that
\begin{align*}
    C_\Phi(\overline{X}) \vdash \pite Y. \Phi(\overline{X}, Y) &  & \text{ and } & & \pite Y. \Phi(\overline{X}, Y) \vdash C_\Phi(\overline{X})
\end{align*}
both hold as claimed.
\end{proof}
\end{proposition}

\begin{body}
Informally, one can summarize Proposition~\ref{prop:general:definability} as follows: any quantifier-free formula that defines $C_\Phi(\overline{X})$ \textit{a fortiori} constitutes a uniform interpolant. Since any two uniform interpolants are equivalent, the result follows. Keeping in mind this connection between $\pite Y. \Phi(\overline{X},Y)$ and the definability of $C_\Phi(\overline{X})$, we can proceed to introduce the concept of auxiliary formula. Analyzing these formulae will serve as our primary strategy when proving the nondefinability of specific regular connectives.
\end{body}

\begin{definition}\label{def:general:auxiliary}
We call the formula $\varphi$ an \textit{auxiliary formula} for the regular connective $C_\Phi(\overline{X})$ if $\pite Y. \Phi(\overline{X},Y) \vdash \Phi(\overline{X},\varphi)$ is derivable.
\end{definition}

\begin{example}\label{example:general:auxiliaryconj}
Take the regular connective $C_\Phi(P,Q)$ defined by the following formula: $(Y \vee \neg Y) \rightarrow (P \wedge Q)$. The formula $P \wedge Q$ constitutes an auxiliary formula for this connective.
\begin{proof}
We observe that the connective $C_\Phi(P,Q)$ is in fact definable in intuitionistic propositional logic by the formula $\neg\neg (P \wedge Q)$. For $\neg\neg (P \wedge Q) \vdash \exists Y. (Y \vee \neg Y) \rightarrow (P \wedge Q)$, notice that one can take $Y$ as $P \wedge Q$ itself; the other direction requires only a straightforward verification. By Proposition~\ref{prop:general:definability}, $\pite Y. (Y \vee \neg Y) \rightarrow (P \wedge Q)$ therefore coincides with $\neg\neg (P \wedge Q)$, and $P \wedge Q$ constitutes an auxiliary formula for the connective.
\end{proof}
\end{example}

\begin{theorem}\label{thm:general:auxiliaryequiv}
If a regular connective $C_\Phi(\overline{X})$ has an auxiliary formula, it is definable in intuitionistic propositional logic. Conversely, if $\pite Y. \Phi(\overline{X}, Y)$ is equivalent to a $\vee$-free formula in intuitionistic propositional logic and $\pite Y. \Phi(\overline{X}, Y) \vdash C_\Phi({\overline{X}})$, then $C_\Phi(\overline{X})$ has an auxiliary formula. Moreover, any $C_\Phi(\overline{X})$ that has some auxiliary formula also has a quantifier-free auxiliary formula.
\begin{proof}
If $C_\Phi(\overline{X})$ has an auxiliary formula $\varphi(\overline{X})$, then $\pite Y. \Phi(\overline{X},Y) \vdash \Phi(\overline{X},\varphi)$ is derivable. Applying the $\exists R$ rule, $\pite Y. \Phi(\overline{X},Y) \vdash C_\Phi(\overline{X})$ follows immediately. But we already know from Theorem~\ref{thm:general:pitts} that $C_\Phi(\overline{X}) \vdash \pite Y. \Phi(\overline{X},Y)$ holds. This gives the definability of $C_\Phi(\overline{X})$.
For the converse, w.l.o.g. assume that $\pite Y. \Phi(\overline{X}, Y)$ is $\vee$-free. Consider a cut-free proof of $\pite Y. \Phi(\overline{X}, Y) \vdashat C_\Phi({\overline{X}})$. We can extract an appropriate $\varphi$ from the proof using the following recursive procedure: climb up the proof tree until the first application of either a $wR$, $\exists R$, $\rightarrow L$ rule. Note that no branching or irreversible rules can occur until that point. Proceed by cases:
\begin{itemize}
    \item $wR$: The $wR$ rule is applied to some $\Gamma' \vdashat C_\Phi(\overline{X})$, yielding $\Gamma' \vdashat$. This means that $\Gamma' \vdashat \Phi(\overline{X},\bot)$ is provable by an application of weakening. And since $\Gamma'$ was obtained from $\pite Y. \Phi(\overline{X},y)$ via reversible rules, we have $\pite Y. \Phi(\overline{X}, Y) \vdash \Phi(\overline{X},\bot)$, and thus $\bot$ is an auxiliary formula for the connective.
    \item $\exists R$: The $\exists R$ rule is applied to some $\Gamma' \vdashat C_\Phi(\overline{X})$, yielding $\Gamma' \vdashat \Phi(\overline{X},\varphi)$ for some quantifier-free formula $\varphi$. Again, by reversibility,  $\pite Y. \Phi(\overline{X}, Y) \vdash \Phi(\overline{X},\varphi)$ and $\varphi$ constitutes an auxiliary formula for the connective $C_\Phi(\overline{X})$.
    \item $\rightarrow L$: The $\rightarrow L$ rule is applied to some $\Gamma', A \rightarrow B \vdashat C_\Phi(\overline{X})$, yielding $\Gamma' \vdashat A$ and $\Gamma', B \vdashat C_\Phi(\overline{X})$. Recurse on the latter to extract a $\varphi$ so that $\Gamma', B \vdashat \Phi(\overline{X}, \varphi)$. Since we have both $\Gamma' \vdashat A$ and $\Gamma', B \vdashat \Phi(\overline{X}, \varphi)$, the $\rightarrow L$ rule gives $\Gamma', A \rightarrow B \vdashat \Phi(\overline{X}, \varphi)$, and since all the rules applied between $\pite Y. \Phi(\overline{X}, Y) \vdashat C_\Phi({\overline{X}})$ and $\Gamma', A \rightarrow B \vdashat C_\Phi(\overline{X})$ were reversible, we get $\pite Y. \Phi(\overline{X}, Y) \vdashat \Phi(\overline{X}, \varphi)$. Thus $\varphi$ is an auxiliary formula for the connective.
\end{itemize}
Finally, we prove that if $C_\Phi(\overline{X})$ has some auxiliary formula $\varphi$, then it in fact has a quantifier-free auxiliary formula. Assume that $\pite Y. \Phi(\overline{X}, Y) \vdashat \Phi(\overline{X}, \varphi)$ has a derivation. By the result of Pitts we have an effectively computable translation $(-)^p$ from the second-order calculus to its quantifier-free fragment which respects derivability and restricts to the identity over quantifier-free formulae. Since $\Phi(\overline{X}, Y)$ itself is quantifier-free by the definition of regular connective, and the left-hand side of $\pite Y. \Phi(\overline{X}, Y) \vdashat \Phi(\overline{X}, \varphi)$ is quantifier-free by Theorem~\ref{thm:general:pitts}, applying the translation yields the derivability of $\pite Y. \Phi(\overline{X}, Y) \vdashat \Phi(\overline{X}, (\varphi)^p)$. Thus $(\varphi)^p$ constitutes a quantifier-free auxiliary formula for $C_\Phi(\overline{X})$.
\end{proof}
\end{theorem}

%%%%%%%%%%%%%%%%%%%%%%%%%%%%%%%%%%%%%%%%%%%%%%%%%%%%%%%%%%%%%%%%%%%%%%%%5

\section{Realizability disjunction}\label{sec:tara}

\begin{definition}\label{def:tara:disjunction}
The regular connective $C_\Phi(X_1,X_2)$ defined by the formula
$$(\neg Y \rightarrow X_1) \wedge (\neg\neg Y \rightarrow X_2)$$
is called \textit{realizability disjunction} and denoted $X_1 \parr X_2$.
\end{definition}

\begin{body}
Realizability disjunction was first defined by Taranovsky~\cite{taranovsky-nonconstructive} under the name \textit{nonconstructive disjunction}. The motivation comes from realizability. Normally, a realizer for the formula $A\vee B$ consists of a pair of objects $(n,r)$: a natural number $n\in\{0,1\}$ and a realizer $r$ which realizes $A$ if $n=0$ and $B$ if $n=1$. So a realizer for a disjunction explicitly indicates one of the disjuncts using the number $n$ and provides a realizer $r$ for the indicated disjunct. Taranovsky suggests a new sort of disjunction-like formula $A \parr B$, with markedly different realizability rules: instead of a pair $(n,r)$, a realizer for $A \parr B$ consists of a pair of realizers $(a,b)$, so that ``$a$ does not realize $A$'' and ``$b$ does not realize $B$'' are not both false. So the realizers for $A \parr B$ have diminished constructive content compared to those for $A \vee B$, as they do not indicate which disjunct was realized.
\end{body}

\begin{body}
Bauer~\cite{bauer-equivalence} showed that in a realizability topos that interprets higher-order logic, one has an equivalence between Taranovsky's realizability-theoretic definition of $\parr$ and the regular connective given in Definition~\ref{def:tara:disjunction}. Thus, in the realizability setting, $\parr$ gives rise to a bona-fide regular second-order definable connective. Bauer asked whether one could find a first-order schema equivalent to $X_1 \parr X_2$.
\end{body}

\begin{body}
Abandoning its origins in realizability, we treat $\parr$ as an ordinary regular connective of second-order propositional logic, given by the formula above (Definition~\ref{def:tara:disjunction}). We apply the method outlined in the introduction to investigate the question of its definability using a quantifier-free schema, and provide a negative answer to Bauer's question in this setting.
\end{body}

\subsection{Elementary properties}

\begin{body}
We begin our analysis with a series of propositions on the properties of $\parr$. This includes establishing its elementary properties, such as commutativity, then confirming that $\parr$ satisfies an analogue of weak excluded middle: $\neg P \parr \neg\neg P$ always holds. We prove a claim characterizing $\parr$ as the strongest among monotonic binary connectives with these properties~(Proposition~\ref{prop:tara:characterization}). All the results presented up to that point were identified by Taranovsky in his initial proposal, albeit without accompanying proofs. In Proposition~\ref{prop:tara:ultimate} we show the implication $(P \parr Q) \rightarrow \neg P \rightarrow Q$. As the Pitts interpolant calculation will reveal, this holds the key to the connective's non-definability (see~\ref{sec:tara:key}). We argue informally, leaving it to the reader to translate the proofs into the appropriate sequent calculus derivation trees.
\end{body}

\begin{proposition}\label{prop:tara:commutative}
The connective $\parr$ is commutative: if $P \parr Q$ holds, so does $Q \parr P$.
\begin{proof}
Assume $P \parr Q$ holds, i.e. we can find some $Y$ such that $\neg Y \rightarrow P$ and $\neg\neg Y \rightarrow Q$ both hold. We construct a $Z$ such that $\neg Z \rightarrow Q$ and $\neg\neg Z \rightarrow P$. We notice immediately that setting $Z$ to $\neg Y$, we get $\neg\neg Y \rightarrow Q$ by assumption. To prove $\neg\neg\neg Y \rightarrow P$, we can invoke triple-negation elimination.
\end{proof}
\end{proposition}

\begin{proposition}\label{prop:tara:monotone}
The connective $\parr$ is monotone: if $Q \rightarrow Q'$ and $P \parr Q$ holds, then so does $P \parr Q'$.
\begin{proof}
Assume $Q \rightarrow Q'$ and $P \parr Q$ both hold. From this, we know that we can find some $Y$ so that $\neg Y \rightarrow P$ and $\neg\neg Y \rightarrow Q$. It suffices to show that $\neg\neg Y \rightarrow Q'$: but that follows from transitivity of implication.
\end{proof}
\end{proposition}

\begin{proposition}\label{prop:tara:bottom}
The propositions $\bot \parr P$ and $P$ are equivalent.
\begin{proof}
For the forward direction, assume that $\bot \parr P$ holds, i.e. we can find some $Y$ so that $\neg Y \rightarrow \bot$ and $\neg\neg Y \rightarrow P$. From $\neg Y \rightarrow \bot$ we know $\neg\neg Y$, and since $\neg\neg Y \rightarrow P$, we can deduce $P$. For the backward direction, assume $P$. Then we have both $\neg P \rightarrow \bot$ and $\neg\neg P \rightarrow P$. Thus, $\bot \parr P$ holds as claimed.
\end{proof}
\end{proposition}

\begin{proposition}\label{prop:tara:wlem}
The connective $\parr$ satisfies an analogue of weak excluded middle: for each proposition $P$, we have $\neg P \parr \neg\neg P$.
\begin{proof}
Follows immediately from the fact that $\neg P \rightarrow \neg P$ and $\neg\neg P \rightarrow \neg\neg P$.
\end{proof}
\end{proposition}

\begin{proposition}\label{prop:tara:characterization}
The $\parr$ connective is strongest among the monotonic (in both arguments) binary connectives that satisfy the analogue of weak excluded middle. In other words, if we have a connective $\oplus$ that satisfies monotonicity and $\neg P \oplus \neg\neg P$, then $P \parr Q \rightarrow P \oplus Q$.
\begin{proof}
Assume that $P \parr Q$ holds. Then we have some $Y$ such that $\neg Y \rightarrow P$ and $\neg\neg Y \rightarrow Q$. Since $\oplus$ satisfies the analogue of weak excluded middle in general, we have that $\neg Y \oplus \neg\neg Y$ for that particular $Y$. Monotonicity, with $\neg Y \rightarrow P$ implies $P \oplus \neg\neg Y$, and with $\neg\neg Y \rightarrow Q$ implies $P \oplus Q$. 
\end{proof}
\end{proposition}

\begin{proposition}\label{prop:tara:ultimate}
The implication $(P \parr Q) \rightarrow \neg P \rightarrow Q$ holds.
\begin{proof}
Assume $P \parr P$ and $\neg P$. We can then find $Y$ so that $\neg Y \rightarrow P$ and $\neg\neg Y \rightarrow Q$ both hold. By assumption $\neg P$ holds, so we have $\neg Y \rightarrow (P \wedge \neg P)$, and consequently $\neg \neg Y$. But since $\neg\neg Y \rightarrow Q$, we also have $Q$. Discharging the assumption $\neg P$, we obtain the desired implication.
\end{proof}
\end{proposition}

\subsection{Non-definability}

\begin{body}
In this section, we apply the method outlined in Section~\ref{sec:introduction} to deal with the question of definability for realizability disjunction. Our investigation culminates in the result that the connective indeed eludes a definition by a quantifier-free schema.
\end{body}

\begin{proposition}\label{prop:tara:calculation}
The formula $\pite Y. (\neg Y \rightarrow X_1) \wedge (\neg\neg Y \rightarrow X_2)$ is equivalent to $$(\neg X_1 \rightarrow X_2) \wedge (\neg X_2 \rightarrow X_1).$$
\begin{proof}
Computation using the \texttt{propquant} tool of F{\'{e}}r{\'{e}}e and Van Gool.
\end{proof}
\end{proposition}

\begin{body}\label{sec:tara:key}
Note that Proposition~\ref{prop:tara:ultimate} provides a quick sanity check on the calculation of Proposition~\ref{prop:tara:calculation}: we already knew that $\neg P \rightarrow Q$ and $\neg Q \rightarrow P$ follow from $P \parr Q$, and the calculation of the Pitts interpolants has now shown that any quantifier-free formula that follows from $P \parr Q$ in fact follows from these two.
\end{body}

\begin{body}
The proof of Proposition~\ref{prop:tara:calculation} concludes the first step in the outlined method. Now, we need to show that, assuming one can define the connective $\parr$ by a quantifier-free schema, it has an auxiliary formula $t$ satisfying certain strong properties. Since we can write the Pitts interpolant in the $\vee$-free fragment, Proposition~\ref{thm:general:auxiliaryequiv} applies.
\end{body}

\begin{proposition}\label{prop:tara:formulae}
The realizability disjunction connective $\parr$ is definable by a quantifier-free schema within intuitionistic second-order propositional logic precisely if we can find an auxiliary formula $t(P,Q)$ satisfying the following:
\begin{enumerate}
    \item $\neg P \rightarrow Q, \neg Q \rightarrow P, \neg t(P,Q) \vdash P$, and
    \item $\neg P \rightarrow Q, \neg Q \rightarrow P, \neg\neg t(P,Q) \vdash Q$.
\end{enumerate}
\begin{proof}
From Proposition~\ref{prop:general:definability} and the calculation of Proposition~\ref{prop:tara:calculation}, we know that, if one can define $\parr$ at all in intuitionistic propositional logic, then in fact one can define it by the schema $(\neg P \rightarrow Q) \wedge (\neg Q \rightarrow P)$. We already know that
$$ \exists Y. (\neg Y \rightarrow P) \wedge (\neg\neg Y \rightarrow Q) \vdash (\neg P \rightarrow Q) \wedge (\neg Q \rightarrow P)$$
from Theorem~\ref{thm:general:pitts}. Thus, $\parr$ is definable precisely if the converse implication,
$$(\neg P \rightarrow Q) \wedge (\neg Q \rightarrow P) \vdash \exists Y. (\neg Y \rightarrow P) \wedge (\neg\neg Y \rightarrow Q)$$
holds as well. According to Theorem~\ref{thm:general:auxiliaryequiv}, this, in turn, happens if
$$(\neg P \rightarrow Q) \wedge (\neg Q \rightarrow P) \vdash (\neg t(P,Q) \rightarrow P) \wedge (\neg\neg t(P,Q) \rightarrow Q)$$
for some auxiliary formula $t(P,Q)$. This immediately gives rise to the two conditions above.
\end{proof}
\end{proposition}

\begin{body}
Proposition~\ref{prop:tara:formulae} completes the second step of the outlined method. In the third step, having identified the auxiliary formula $t(P,Q)$ and its two properties, we must determine whether its presence allows us to deduce a formula that is not intuitionistically valid. As ever with such tasks, the devil is in the details of the syntactic reasoning.
\end{body}

\begin{theorem}\label{thm:tara:mainresult}
Assume that we can find an intuitionistic propositional formula $t(P,Q)$ in two variables so that
\begin{enumerate}
    \item $\neg P \rightarrow Q, \neg Q \rightarrow P, \neg t(P,Q) \vdashl P$, and
    \item $\neg P \rightarrow Q, \neg Q \rightarrow P, \neg\neg t(P,Q) \vdashl Q$.
\end{enumerate}
in a super-intuitionistic logic $\mathcal{L}$. Then $\mathcal{L}$ coincides with classical propositional logic.
\end{theorem}

\subsection{Proof of Theorem~\ref{thm:tara:mainresult}}

\begin{body}
We aim to show that the auxiliary term $t(P,Q)$ obeys rules analogous to those of implication, including Modus Ponens (Proposition~\ref{prop:tara:implication}) and also properties similar to the classical rules for negated conditionals, such as $\neg\neg t(\top, P) \vdashl P$. This forces implication itself to obey a similar rule in $\mathcal{L}$, and double-negation elimination follows as a consequence (Theorem~\ref{thm:tara:dnegelim}).
\end{body}
\pagebreak

\begin{proposition}\label{prop:tara:topreducts}
We can infer all of
\begin{enumerate}
    \item $t(P,Q), Q \vdashl t(P,\top)$
    \item $t(P,Q), P \vdashl t(\top,Q)$
    \item $t(P,\top), Q \vdashl t(P,Q)$
    \item $t(\top,Q), P \vdashl t(P,Q)$
\end{enumerate}
for any formulae $P,Q$.
\begin{proof}
\begin{align*}
 & \textbf{No.} & & \textbf{Claim} & & \textbf{Justification} \\
 & 1 & &
 Q \rightarrow \top, \top \rightarrow Q, t(P,Q) \vdashl t(P, \top) & &
 \text{Lemma~\ref{lemma:general:extensionality}} \\
 & 2 & &
 t(P,Q), Q \vdash Q \rightarrow \top & &
 \text{ } \\
 & 3 & &
 t(P,Q), Q \vdash \top \rightarrow Q & &
 \text{ } \\
 & 4 & &
 t(P,Q), Q \vdashl t(P, \top) & &
 \text{cuts on 1,2,3} \\
\end{align*}
takes care of $t(P,Q), Q \vdashl t(P, \top)$. Nearly identical arguments handle all the other cases as well.
\end{proof}
\end{proposition}

\begin{proposition}\label{prop:tara:negtopreducts}
We can infer $\neg t(P,\top) \vdashl P$ and $\neg\neg t(\top,Q) \vdashl Q$ for any formulae $P,Q$.
\begin{proof}
Since $\neg P \rightarrow \top, \neg \top \rightarrow P, \neg t(P,\top) \vdashl P$ constitutes an instance of the first defining property of $t(P,Q)$, cutting against $\vdash \neg P \rightarrow \top$ and $\vdash \neg \top \rightarrow P$ gives $\neg t(P,\top) \vdashl P$. A similar argument using the second defining property gives $\neg\neg t(\top,Q) \vdashl Q$.
\end{proof}
\end{proposition}

\begin{proposition}\label{prop:tara:implication}
We can infer $t(P,Q), P \vdashl Q$ for any formulae $P,Q$.
\begin{proof}
Use Propositions~\ref{prop:tara:topreducts} and \ref{prop:tara:negtopreducts} as follows.
\begin{align*}
 & \textbf{No.} & & \textbf{Claim} & & \textbf{Justification} \\
 & 1 & &
 t(P,Q), P \vdashl t(\top, Q) & &
 \text{Proposition~\ref{prop:tara:topreducts}} \\
 & 2 & &
 t(\top, Q) \vdash \neg\neg t(\top, Q) & &
 \text{ } \\
 & 3 & &
 \neg\neg t(\top,Q) \vdashl Q & &
 \text{Proposition~\ref{prop:tara:negtopreducts}} \\
 & 4 & &
 t(P, Q), P \vdashl Q & &
 \text{cuts on 1,2,3} \\
\end{align*}
~
\end{proof}
\end{proposition}

\begin{proposition}\label{prop:tara:topswitch}
We can infer $t(P, \top) \wedge P \vdashl t(\top,P)$ and $t(\top,P) \vdashl P \wedge t(P,\top)$ for any formula $P$. In other words, the formulae $t(P, \top) \wedge P$ and $t(\top, P)$ are equivalent, so by Lemma~\ref{lemma:general:extensionality}, one can replace the other in any formula while preserving derivability.
\begin{proof}
One direction is immediate from Proposition~\ref{prop:tara:topreducts}. For the other direction, we can argue as follows:

\begin{align*}
 & \textbf{No.} & & \textbf{Claim} & & \textbf{Justification} \\
 & 1 & &
 t(\top, P) \vdash \neg\neg t(\top, P) & &
 \text{ } \\
 & 2 & &
 \neg\neg t(\top, P) \vdashl P & &
 \text{Proposition~\ref{prop:tara:negtopreducts}} \\
 & 3 & &
 t(\top, P) \vdash P & &
 \text{cut on 1,2} \\
 & 4 & &
 t(\top, P), P \vdash P & &
 \text{ } \\
 & 5 & &
 t(\top, P), P \vdashl t(P,P) & &
 \text{Proposition~\ref{prop:tara:topreducts}} \\
 & 6 & &
 t(\top, P), P \vdashl t(P,P) \wedge P & &
 \text{$\wedge R$ on 4,5} \\
 & 7 & &
 t(P,P), P \vdashl t(P, \top) & &
 \text{Proposition~\ref{prop:tara:topreducts}} \\
 & 8 & &
 t(P,P) \wedge P \vdashl t(P, \top) & &
 \text{$\wedge L$ on 7} \\
 & 9 & &
 t(\top, P), P \vdashl t(P, \top) & &
 \text{cut on 6,8} \\
 & 10 & &
 t(\top, P) \vdashl t(P, \top) & &
 \text{cut on 3,9} \\
 & 11 & &
 t(\top, P) \vdashl t(P, \top) \wedge P & &
 \text{$\wedge R$ on 3,10} \\
\end{align*}
\end{proof}
\end{proposition}

\begin{proposition}\label{prop:tara:fulcrum}
We can infer $\neg t(P,\neg t(\top,P)) \vdashl P$ and $\neg t(P,P \rightarrow \neg t(P,\top)) \vdashl P$ for any formula $P$.
\begin{proof}
The first claim can be derived as follows:
\begin{align*}
 & \textbf{No.} & & \textbf{Claim} & & \textbf{Justification} \\
 & 1 & &
 \neg P \rightarrow \neg t(\top, P), \neg \neg t(\top, P) \rightarrow P, \neg t(P,\neg t(\top, P)) \vdashl P & &
 \text{def. prop. 2} \\
 & 2 & &
 \neg\neg t(\top,P) \vdashl P & &
 \text{Proposition~\ref{prop:tara:negtopreducts}} \\
 & 3 & &
 \vdashl \neg\neg t(\top,P) \rightarrow P & &
 \text{$\rightarrow R$ on 2} \\
 & 4 & &
 \neg\neg t(\top,P) \rightarrow P \vdash \neg P \rightarrow \neg t(\top,P) & &
 \text{ } \\
 & 5 & &
 \vdashl \neg P \rightarrow \neg t(\top,P) & &
 \text{cut on 3,4} \\
  & 6 & &
 \neg t(P,\neg t(\top, P)) \vdashl P & &
 \text{cuts on 1,3,5} \\
\end{align*}
For the second claim, observe that $P \rightarrow \neg t(P,\top)$ is intuitionistically equivalent to $\neg (P \wedge t(P,\top))$, and therefore
$\neg t(P, P \rightarrow \neg t(P,\top)) \vdash \neg t(P, \neg (P \wedge t(P,\top)))$. Similarly, we already know from Proposition~\ref{prop:tara:topswitch} that $P \wedge t(P,\top)$ is $\mathcal{L}$-equivalent to $t(\top, P)$. Thus $\neg t(P, \neg (P \wedge t(P,\top))) \vdashl \neg t(P, \neg t(\top, P))$. At this point, applying the first claim gives $\neg t(P,P \rightarrow \neg t(P,\top)) \vdashl P$ as desired.
\end{proof}
\end{proposition}

\begin{theorem}\label{thm:tara:dnegelim}
We can infer $\neg\neg P \vdashl P$ for any formula $P$.
\begin{proof}
As shown in Proposition~\ref{prop:tara:fulcrum}, we have $\neg t(P,P\rightarrow \neg t(P,\top)) \vdash P$. We show that $\neg \neg P \vdashl \neg t(P,P\rightarrow \neg t(P,\top))$, and consequently $\neg\neg P \vdashl P$ as claimed. Thanks to contrapositives, it suffices to establish $t(P,P\rightarrow \neg t(P,\top)) \vdashl \neg P$ as follows:
\begin{align*}
 & \textbf{No.} & & \textbf{Claim} & & \textbf{Justification} \\
 & 1 & &
 t(P,P\rightarrow \neg t(P,\top)), P \vdashl P \rightarrow \neg t(P,\top) & &
 \text{Proposition~\ref{prop:tara:implication}} \\
 & 2 & &
 t(P,P\rightarrow \neg t(P,\top)), P \vdash P & &
 \text{ } \\
 & 3 & &
 t(P,P\rightarrow \neg t(P,\top)), P \vdashl P \wedge (P \rightarrow \neg t(P,\top)) & &
 \text{$\wedge R$ on 1,2} \\
 & 4 & &
 P \wedge (P \rightarrow \neg t(P,\top)) \vdash \neg t(P,\top) & &
 \text{ } \\
 & 5 & &
 t(P,P\rightarrow \neg t(P,\top)), P \vdashl \neg t(P,\top) & &
 \text{cut on 3,4} \\
 & 6 & &
 t(P,P\rightarrow \neg t(P,\top)) \vdashl P \rightarrow \neg t(P,\top) & &
 \text{$\rightarrow R$ on 1e} \\
 & 7 & &
 t(P,P\rightarrow \neg t(P,\top)), P \rightarrow \neg t(P,\top) \vdashl t(P,\top) & &
 \text{Proposition~\ref{prop:tara:topreducts}} \\
 & 8 & &
 t(P,P\rightarrow \neg t(P,\top)) \vdashl t(P,\top) & &
 \text{cut on 6,7} \\
 & 9 & &
 t(P,P\rightarrow \neg t(P,\top)), P \vdashl t(P,\top) & &
 \text{$wL$ on 8} \\
 & 10 & &
 t(P,P\rightarrow \neg t(P,\top)), P \vdashl \bot & &
 \text{from 6,9} \\
 & 11 & &
 t(P,P\rightarrow \neg t(P,\top)) \vdashl \neg P & &
 \text{$\neg R$ on 10} \\
\end{align*}
Since $\neg\neg P \vdashl P$, we conclude that $\mathcal{L}$ coincides with classical propositional logic.
\end{proof}
\end{theorem}

\begin{body}
From Proposition~\ref{prop:tara:formulae} and Theorem~\ref{thm:tara:dnegelim} we immediately get the non-definability of the $\parr$ connective: if intuitionistic propositional logic could define it, then by Theorem~\ref{thm:general:auxiliaryequiv} its auxiliary term would satisfy the defining properties above, allowing us to construct a proof of $\neg\neg P \vdash P$ for any $P$ inside intuitionistic propositional logic itself.
\end{body}

\subsection{Results in type theory}\label{sec:typetheory}

\begin{body}
The \textit{propositions-as-some-types} perspective came to prominence with the advent of univalent type theories such as Homotopy Type Theory, and identifies propositions as corresponding to "types with at most one inhabitant" (the \textit{h-propositions} in homotopy jargon). This view stands in contrast to the propositions-as-types paradigm that considers every type a proposition, and its inhabiting terms as proofs. We obtain our type-theoretic results in the propositions-as-some-types paradigm. We have a type of all propositions ($\mathtt{Prop}$), and our logical connectives and formulae, including the analogue of the auxiliary formula $t$ above, map propositions to other propositions. This retains compatibility with the choices made in The HoTT Book~\cite{hottbook}, so the reader not well-versed in type theory can use it as an introduction or reference for our type-theoretic results.
\end{body}

\begin{body}
Proof assistants (also known as interactive theorem provers) are computer software that aid mathematicians and computer scientists in defining formal mathematical theories, constructing proofs, and checking their correctness. Agda~\cite{norell-agda} is one of the many proof assistants built upon dependent type theory. Agda works as an interactive, as opposed to an automated theorem prover: it does not generate the proof by itself, but verifies proof scripts that have been encoded into its programming language by the mathematician.
\end{body}

\begin{body}
One can \textit{execute} pure proof-theoretic arguments, such as the proof of Theorem~\ref{thm:tara:mainresult} given above directly inside the logic of a proof assistant. Formalizing the sequent calculus within Agda, ensuring the fidelity of cut-elimination theorems, etc. would involve intense effort. We do not formalize the mechanics of sequent calculus within Agda, and we do not prove the existence of a deduction; instead, we directly deduce $\neg\neg P \rightarrow P$ \textit{within} Agda's native logic, under the assumption that terms inhabiting the types corresponding to the two assumptions exist. We use a very minimal setting, plain Agda with Escardó's $\mathtt{Prop}$ type implementation, as our basis, but the same results can be replayed in "Book HoTT", Cubical Type Theory and nearly all other systems following the propositions-as-some-types paradigm.  Translating our theorem into Agda's internal logic effectively gets us a new theorem with minimal added effort: while the result stated about sequent calculi above and the type-theoretic one formally verified by the proof assistant are closely related, they are ultimately different statements. The full formalization can be found in the Git repository hosted at
\begin{center}
\url{https://github.com/zaklogician/proof-theoretic-methods}.
\end{center}
We state the relevant main result as Proposition~\ref{prop:tara:univalent}: for the corresponding proof, see the file \texttt{RDisjunction.agda} of the aforementioned repository.
\end{body}

\begin{proposition}[in Agda]\label{prop:tara:univalent}
Assume that we are given some $t : \mathtt{Prop} \rightarrow \mathtt{Prop} \rightarrow \mathtt{Prop}$ and inhabitants of the following types:
\begin{itemize}
    \item $\Pi (P\:Q : \mathtt{Prop}). ((P \Rightarrow \bot) \Rightarrow Q) \Rightarrow ((Q \Rightarrow \bot) \Rightarrow P) \Rightarrow (t\:P\:Q \Rightarrow \bot) \Rightarrow P$, and
    \item $\Pi (P\:Q : \mathtt{Prop}). ((P \Rightarrow \bot) \Rightarrow Q) \Rightarrow ((Q \Rightarrow \bot) \Rightarrow P) \Rightarrow ((t\:P\:Q \Rightarrow \bot) \Rightarrow \bot) \Rightarrow Q$.
\end{itemize}
Then we can construct an inhabitant of the type $\Pi P: \mathtt{Prop}. ((P \Rightarrow \bot) \Rightarrow \bot) \Rightarrow P$.
\end{proposition}

%%%%%%%%%%%%%%%%%%%%%%%%%%%%%%%%%%%%%%%%%%%%%%%%%%%%%%%%%%%%%%%%%%%%%%%%

\section{Kreisel's star connective}\label{sec:kreisel}

\begin{definition}\label{def:kreisel:star}
The regular connective $C_\Phi(X)$ defined by the formula $X \leftrightarrow (\neg Y \vee \neg\neg Y)$ is called \textit{Kreisel's star connective}, and denoted $\ast(X)$.
\end{definition}

\begin{body}
The star connective of Definition~\ref{def:kreisel:star} was introduced by Kreisel~\cite{kreisel-monadic}, who made use of topological semantics to prove that $\ast(X)$ is not definable in intuitionistic propositional logic (in contrast to its counterpart $X \leftrightarrow (Y \vee \neg Y)$, which one can define using the formula $\neg\neg X$). This was also the first result to imply that the $\exists$ quantifier itself is not definable in quantifier-free intuitionistic propositional logic. Troelstra~\cite{troelstra-kreiselstar}, building upon Kreisel's foundation, gives examples of topological spaces in which one cannot define \textit{any} of the regular connectives given by $X \leftrightarrow \psi(Y)$ where $\psi(Y)$ denotes a quantifier-free formula in one variable $Y$ such that $\neg Y \vee \neg\neg Y \vdash \psi(Y)$. It is this more general nondefinability result which we re-establish using proof-theoretic methods in Theorem~\ref{thm:kreisel:troelstra}, allowing us to strengthen the conclusion so that it applies to a wide range of super-intuitionistic logics (Corollaries~\ref{cor:kreisel:extensions}, \ref{cor:kreisel:disjunction}).
\end{body}

\begin{body}\label{body:topcomments}
Troelstra~\cite{troelstra-kreiselstar} notes the nuanced dependence of the $\ast(P)$ connective's definability on the particular topological space under consideration. For instance, in the open interval $(0,1)$ and in Cantor space, $\ast(P)$ coincides with $\neg\neg P$: this is in strict contrast to the quantifier-free setting, where such \textit{dense-in-itself} spaces suffice to provide a complete semantics for the logic. A later contribution by Po\l acik~\cite{polacik-pitts} greatly elucidates how dense-in-itself metric spaces work for second-order logic: in particular, assuming a classical metatheory, nullary regular connectives are always definable in them.
\end{body}

\subsection{General non-definability}

\begin{body}
We apply the same strategy as before. F{\'{e}}r{\'{e}}e-Van Gool's \texttt{propquant} tool calculates $\pite Y. P \leftrightarrow (\neg Y \vee \neg\neg Y)$ as $\neg\neg P$ immediately. Since the interpolant only has one variable, the Rieger-Nishimura lattice can also be used to manually double-check its value.
\end{body}

\begin{proposition}\label{prop:kreisel:riegerlower}
Consider a quantifier-free formula in one variable $\psi(Y)$ so that $\neg Y \vee \neg\neg Y \vdash \psi(Y)$. All of the following hold:
\begin{enumerate}
    \item $\vdash \neg\neg \psi(Y)$,
    \item $Y \vdash \psi(Y)$,
    \item $\vdash \psi(\psi(Y))$, and
    \item $\psi(Y) \rightarrow Y \vdash Y$.
\end{enumerate}
\begin{proof}
From $\neg Y \vee \neg\neg Y \vdash \psi(Y)$ we know that $\psi(Y)$ is a classical tautology, so we have $\vdash \neg\neg \psi(Y)$ by Glivenko's theorem. Cutting $Y \vdash \neg Y \vee \neg\neg Y$ against $\neg Y \vee \neg\neg Y \vdash \psi(Y)$, we have $Y \vdash \psi(Y)$. We obtain $\psi(\psi(Y))$ as follows:
\begin{align*}
 & \textbf{No.} & & \textbf{Claim} & & \textbf{Justification} \\
 & 1 & &
 \neg \psi(Y) \vee \neg\neg \psi(Y) \vdash \psi(\psi(Y)) & &
 \text{assm.} \\
 & 2 & &
 \vdash \neg\neg \psi(Y) & &
 \text{ } \\
 & 3 & &
 \vdash \neg \psi(Y) \vee \neg\neg \psi(Y) & &
 \text{$\vee R_2$ on 2} \\
 & 4 & &
 \vdash \psi(\psi(Y)) & &
 \text{cut on 1,3} \\
\end{align*}
We know that $Y \rightarrow \psi(Y), \psi(Y) \rightarrow Y, \psi(\psi(Y)) \vdash \psi(Y)$. Cut this against $\vdash Y \rightarrow \psi(Y)$ and $\vdash \psi(\psi(Y))$ to deduce $\psi(Y) \rightarrow Y \vdash \psi(Y)$. Finally, cut against the tautological $\psi(Y) \rightarrow Y, \psi(Y) \vdash Y$ to obtain the result.
\end{proof}
\end{proposition}

\begin{proposition}\label{prop:kreisel:riegerhigher}
Consider a quantifier-free formula in one variable $\psi(Y)$ which satisfies $\neg Y \vee \neg\neg Y \vdash \psi(Y)$. Then $\pite Y. X \leftrightarrow \psi(Y)$ coincides with $\neg\neg X$ in intuitionistic propositional logic.
\end{proposition}

\begin{theorem}[Troelstra]\label{thm:kreisel:troelstra}
Take a quantifier-free formula $\psi(Y)$ in one variable which satisfies $\neg Y \vee \neg\neg Y \vdash \psi(Y)$. One cannot define the regular connective given by the formula $X \leftrightarrow \psi(Y)$ in intuitionistic propositional logic.
\end{theorem}

\subsection{Proof of Theorem~\ref{thm:kreisel:troelstra}}

\begin{body}
By Proposition~\ref{prop:kreisel:riegerhigher} and Theorem~\ref{thm:general:auxiliaryequiv}, if we could define the regular connective given by the formula $P \leftrightarrow \psi(Y)$, we could find some auxiliary term $t(P)$ satisfying $\neg\neg P \vdash P \leftrightarrow \psi(t(P))$.
\end{body}

\begin{body}
In the rest of this section we work in a logic $\mathcal{L}$ that contains a term $t(Y)$ with the following defining properties:
\begin{enumerate}
    \item $P \vdashl \psi(t(P))$,
    \item $\neg \neg P, \psi(t(P)) \vdashl P$
\end{enumerate}
To aid in parenthesis management we introduce the abbreviations $\psi t(P)$ and $t\psi(P)$ standing for $\psi(t(P))$ and $t(\psi(P))$ respectively.
\end{body}

\begin{body}
The idea of the proof is to derive $\vdashl \psi(t(P))$: combined with the second defining property, double-negation elimination follows as a consequence (Theorem~\ref{thm:kreisel:mainresult}).
\end{body}

\begin{lemma}[Trivium]\label{lemma:kreisel:psitpsi}
All of the following hold:
\begin{enumerate}
    \item $\psi t \psi (P) \vdashl \psi(P)$, 
    \item $\psi(P) \vdashl \psi t \psi (P)$,
    \item $\psi(P) \rightarrow t \psi (P) \vdashl t \psi (P)$.
\end{enumerate}
\begin{proof}The first claim is immediate from the defining property $\neg \neg P, \psi(t(P)) \vdashl P$ and Proposition~\ref{prop:kreisel:riegerlower}. The second claim comes about as an instance of the defining property $P \vdashl \psi t P$. The third claim requires more work:
\begin{align*}
 & \textbf{No.} & & \textbf{Claim} & & \textbf{Justification} \\
 & 1 & &
 \psi t \psi (P) \rightarrow t \psi (P) \vdash t \psi (P) & &
 \text{Proposition~\ref{prop:kreisel:riegerlower}} \\
 & 2 & &
 \psi(P) \vdashl \psi t \psi (P) & &
 \text{Trivium 2} \\
 & 3 & &
 \psi t \psi (P) \rightarrow t \psi (P) \vdashl  \psi(P) \rightarrow t \psi (P) & &
 \text{contrapos. of 2} \\
 & 4 & &
 \psi(P) \rightarrow t \psi (P) \vdashl t \psi(P) & &
 \text{cut on 1,3} \\
\end{align*}
\end{proof}
\end{lemma}

\begin{lemma}[Quadrivium]\label{lemma:kreisel:fourtough}
All of the following hold:
\begin{enumerate} 
\item $t \psi (P) \vdashl t(\top)$,
\item $t(\top) \rightarrow t \psi (P) \vdashl \psi (P)$,
\item $t \psi (P) \vdashl \psi (P)$, and
\item $t(\top) \rightarrow \psi (P) \vdashl \psi (P)$.
\end{enumerate}
\begin{proof}
We know that $Q, t \psi (Q) \vdash t(\top)$ holds for any $Q$ by Lemma~\ref{lemma:general:extensionality} and Proposition~\ref{prop:kreisel:riegerlower}. We will use a substitution instance, replacing $Q$ with $t\psi(P)$ to conclude the first claim.
\begin{align*}
 & \textbf{No.} & & \textbf{Claim} & & \textbf{Justification} \\
 & 1 & &
 t \psi (P), t \psi t \psi (P) \vdashl t(\top) & &
 \text{ } \\
 & 2 & &
 \psi (P) \vdashl \psi t \psi(P) & &
 \text{Lemma~\ref{lemma:kreisel:psitpsi}} \\
 & 3 & &
 \psi t \psi(P) \vdashl \psi(P) & &
 \text{Lemma~\ref{lemma:kreisel:psitpsi}} \\
 & 4 & &
 \vdashl \psi(P) \rightarrow \psi t \psi(P) & &
 \text{$\rightarrow R$ on 2} \\
 & 5 & &
 \vdashl \psi t \psi(P) \rightarrow \psi(P) & &
 \text{$\rightarrow R$ on 3} \\
 & 6 & &
 \psi(P) \rightarrow \psi t \psi (P), \psi t \psi (P) \rightarrow \psi(P), t \psi(P) \vdashl  t \psi t \psi(P) & &
 \text{Lemma~\ref{lemma:general:extensionality}} \\
 & 7 & &
 t \psi(P) \vdashl  t \psi t \psi(P) & &
 \text{cuts on 4,5,6} \\
 & 8 & &
 t \psi (P) \vdashl t(\top)  & &
 \text{cut on 1,7} \\
\end{align*}
The second claim is crucial, so we provide ample detail:
\begin{align*}
 & \textbf{No.} & & \textbf{Claim} & & \textbf{Justification} \\
 & 1 & &
 t(\top) \rightarrow t\psi(P), t\psi(P) \rightarrow t(\top), \psi t (\top) \vdash \psi t \psi (P) & &
 \text{Lemma~\ref{lemma:general:extensionality}} \\
 & 2 & &
 t \psi (P) \vdashl t(\top) & &
 \text{claim 1} \\
 & 3 & &
 \vdashl t \psi (P) \rightarrow t(\top) & &
 \text{$\rightarrow R$ on 2} \\
 & 4 & &
 t(\top) \rightarrow t\psi(P), \psi t (\top) \vdashl \psi t \psi (P)& &
 \text{cut on 1,3} \\
 & 5 & &
 \top \vdashl \psi t (\top) & &
 \text{def. prop. 1}\\
 & 6 & &
 \vdash \top & &
 \text{ }\\
 & 7 & &
 \vdashl \psi t (\top) & &
 \text{cut on 5,6}\\
 & 8 & &
 t(\top) \rightarrow t\psi(P) \vdashl \psi t \psi (P)& &
 \text{cut on 4,7} \\
 & 9 & &
 \psi t \psi (P) \vdashl \psi(P) & &
 \text{Lemma~\ref{lemma:kreisel:psitpsi}} \\
 & 10 & &
 t(\top) \rightarrow t\psi(P) \vdashl \psi (P)& &
 \text{cut on 8,9} \\
\end{align*}
The third claim is immediate by cutting $t\psi(P) \vdash t(\top) \rightarrow t\psi(P)$ with the second claim. We prove the fourth and final claim as follows:
\begin{align*}
 & \textbf{No.} & & \textbf{Claim} & & \textbf{Justification} \\
 & 1 & &
 \psi(P), t(\top) \vdash t \psi (P) & &
 \text{Lemma~\ref{lemma:general:extensionality}} \\
 & 2 & &
 \psi(P) \wedge t(\top) \vdash t \psi(P) & &
 \text{$\wedge L$ on 1} \\
 & 3 & &
 t(\top) \rightarrow (\psi(P) \wedge t(\top)) \vdash t(\top) \rightarrow t \psi(P)& &
 \text{monotonicity on 2} \\
 & 4 & &
 t(\top) \rightarrow \psi(P) \vdash t(\top) \rightarrow (\psi(P) \wedge t(\top)) & &
 \text{ } \\
 & 5 & &
 t(\top) \rightarrow \psi(P) \vdash t(\top) \rightarrow t \psi(P)& &
 \text{cut on 3,4}\\
 & 6 & &
 t(\top) \rightarrow t \psi(P) \vdashl \psi(P) & &
 \text{claim 2}\\
 & 7 & &
 t(\top) \rightarrow \psi(P) \vdashl \psi(P) & &
 \text{cut on 5,6}\\
\end{align*}
This proves all four claims.
\end{proof}
\end{lemma}

\begin{theorem}\label{thm:kreisel:mainresult}
We can infer $\neg\neg P \vdashl P$ for any formula $P$.
\begin{proof}
Both $t \psi t(P) \rightarrow \psi t(P) \vdashl \psi t(P)$ and $\vdashl t \psi t(P) \rightarrow \psi t(P)$ follow from Lemma~\ref{lemma:kreisel:fourtough}, and therefore $\vdashl \psi t(P)$ is derivable. Combining this with the second defining property of $t$ gives double-negation elimination.
\end{proof}
\end{theorem}

\begin{body}
Theorem~\ref{thm:kreisel:troelstra} follows immediately from Theorem~\ref{thm:kreisel:mainresult}: if intuitionistic propositional logic could define the regular connective given by the formula $P \leftrightarrow \psi(Y)$, then by Theorem~\ref{thm:general:auxiliaryequiv} its auxiliary term would satisfy the defining properties above, and hence one could construct a proof of $\neg\neg P \vdash P$ for any $P$. 
\end{body}

\subsection{Consequences}
\begin{body}
The structural proof of Theorem~\ref{thm:kreisel:troelstra} shows that $\neg\neg X$ can essentially never define $\psi(X)$. This observation allows us to prove a strengthened version of Troelstra's original result, which also yields information about super-intuitionistic logics that \textit{can} define these regular connectives.
\end{body}

\begin{corollary}\label{cor:kreisel:extensions}
Consider a quantifier-free formula in one variable $\psi(Y)$ that satisfies $\neg Y \vee \neg\neg Y \vdash \psi(Y)$. Any intermediate logic $\mathcal{L}$ which defines the regular connective given by the formula $X \leftrightarrow \psi(Y)$ proves $\vdashl \psi(P)$.
\begin{proof}
Take such a logic $\mathcal{L}$. We know that $\pite Y. X \leftrightarrow \psi(Y)$ is $\neg\neg X$. Consequently, the unary formula of $\mathcal{L}$ which defines the connective must coincide with either $\neg\neg X$ or $X$. The case where it coincides with the former reduces to one where Theorem~\ref{thm:kreisel:mainresult} applies. Otherwise, $\vdashl X \leftrightarrow \exists Y. X \leftrightarrow \psi(Y)$. Then $\vdash \exists Y. \psi(P) \leftrightarrow \psi(Y)$, so $\vdashl \psi(P)$.
\end{proof}
\end{corollary}

\begin{corollary}\label{cor:kreisel:disjunction}
No intermediate logic that enjoys the disjunction property defines any of the following connectives:
\begin{enumerate}
\item Kreisel's star connective,
\item the regular connective given by $X \leftrightarrow (\neg \neg Y \vee (\neg\neg Y \rightarrow Y))$,
\item the regular connective given by $X \leftrightarrow ((\neg\neg Y \rightarrow Y) \rightarrow (Y \vee \neg Y))$.
\end{enumerate}
\begin{proof}
Apply Corollary~\ref{cor:kreisel:extensions}, and use the fact that the base logics in question already lack the disjunction property.
\end{proof}
\end{corollary}

\begin{body}
Keep in mind that Corollary~\ref{cor:kreisel:disjunction} does not extend arbitrarily high up in the Rieger-Nishimura lattice: Scott's logic famously has the disjunction property.
\end{body}

\begin{body}
As with our proof of Theorem~\ref{thm:tara:mainresult} in the analysis of the realizability disjunction connective $\parr$, the proof of Theorem~\ref{thm:kreisel:mainresult} also applies internally in a univalent type theories. For the development, see the file \texttt{KreiselStar.agda} in the repository; here we only state the main result as Proposition~\ref{prop:kreisel:univalent}.
\end{body}

\begin{proposition}\label{prop:kreisel:univalent}
Univalent type theory proves that whenever we have $t : \mathrm{Prop} \rightarrow \mathrm{Prop}$, $\psi : \mathrm{Prop} \rightarrow \mathrm{Prop}$ so that the types
\begin{enumerate}
  \item $\Pi P: \mathrm{Prop}. P \Rightarrow \psi (t P)$
  \item $\Pi P: \mathrm{Prop}. ((P \Rightarrow \bot) \Rightarrow \bot) \Rightarrow \psi (t P) \Rightarrow P$
  \item $\Pi P: \mathrm{Prop}. (P \Rightarrow \bot) \vee ((P \Rightarrow \bot) \Rightarrow \bot) \Rightarrow \psi P$
\end{enumerate}
have inhabitants, then so has the type $\Pi P: \mathrm{Prop}. ((P \Rightarrow \bot) \Rightarrow \bot) \Rightarrow P$.
\end{proposition}

%%%%%%%%%%%%%%%%%%%%%%%%%%%%%%%%%%%%%%%%%%%%%%%%%%%%%%%%%%%%%%%%%%%%%%%%

\section{Po\l acik's connective}\label{sec:polacik}

\begin{definition}\label{def:polacik:connective}
The regular connective $C_\Phi(X)$ defined by the formula $$(X \rightarrow (Y \vee \neg Y)) \rightarrow X$$ is called \textit{Po\l acik's connective} and denoted $\bullet(X)$.
\end{definition}

\begin{body}
Po\l acik~\cite{polacik-pitts} introduced the connective of Definition~\ref{def:polacik:connective} as an example when studying the relationship between Pitts quantifiers and topological quantification, noting that its undefinability follows by considering a particular Kripke model, yet cannot be shown via topological semantics in any dense-in-itself metric space. The connective saw further use in Zdanowski's work~\cite{zdanowski-existential} characterizing the expressive power of the universal-free fragment of second-order propositional logic. 
\end{body}

\begin{body}
In this brief section, we demonstrate that the method described in Section~\ref{sec:introduction} gives a very short, self-contained proof-theoretic argument showing the non-definability of Po\l acik's connective by a quantifier-free schema.
\end{body}

\begin{proposition}\label{prop:polacik:interpolant}
The formula $\pite Y. (P \rightarrow (Y \vee \neg Y)) \rightarrow P$ coincides with $\neg \neg P$.
\begin{proof}
Via the \texttt{propquant} tool of F{\'{e}}r{\'{e}}e and Van Gool.
\end{proof}
\end{proposition}

\begin{theorem}\label{thm:polacik:undefinable}
One cannot define Po\l acik's connective in intuitionistic propositional logic.
\end{theorem}

\subsection{Proof of Theorem~\ref{thm:polacik:undefinable}}

\begin{body}\label{body:polacikterm}
The proof proceeds the same way as the analogous arguments did for realizability disjunction and Kreisel's star connective, although with substantially reduced complexity compared to the others. By the calculation of Proposition~\ref{prop:polacik:interpolant} and Theorem~\ref{thm:general:auxiliaryequiv}, we know that if we could define the connective $\bullet(P)$, we could find some auxiliary term $t(P)$ satisfying the defining property $$\neg\neg P, P \rightarrow (t(P) \vee \neg t(P)) \vdash P.$$
\end{body}

\begin{body}
The idea is to show that $\vdashl t(P) \vee \neg t(P)$ holds for the auxiliary term described in \ref{body:polacikterm}. Together with the defining property, this yields double-negation elimination for $\mathcal{L}$. 
\end{body}

\begin{body}
We introduce the abbreviation $f(P)$ to stand for $t(P) \vee \neg t(P)$ and further abbreviate $f(f(P))$ as $ff(P)$ and $f(f(f(P)))$ as $fff(P)$ to help with parenthesis management.
\end{body}

\begin{lemma}[Three-in-one]\label{lemma:polacik:three}
We can infer all of the following:
\begin{enumerate}
    \item $f(P), ff(P) \vdashl fff(P)$,
    \item $ff(P) \rightarrow fff(P) \vdashl ff(P)$,
    \item $\vdashl f(P)$.
\end{enumerate}
\begin{proof}
We get the first claim as follows:
\begin{align*}
 & \textbf{No.} & & \textbf{Claim} & & \textbf{Justification} \\
 & 1 & &
 f(P) \rightarrow ff(P), ff(P) \rightarrow f(P), ff(P) \vdashl fff(P) & &
 \text{Lemma~\ref{lemma:general:extensionality}} \\
 & 2 & &
 ff(P) \vdash f(P) \rightarrow ff(P) & &
 \text{ } \\
 & 3 & &
 f(P) \vdash ff(P) \rightarrow f(P) & &
 \text{ } \\
 & 4 & &
 f(P), ff(P) \vdashl fff(P) & &
 \text{cuts on 1,2,3} \\
\end{align*}
Have the second claim by
\begin{align*}
 & \textbf{No.} & & \textbf{Claim} & & \textbf{Justification} \\
 & 1 & &
 \neg\neg ff(P), ff(P) \rightarrow fff(P) \vdashl ff(P) & &
 \text{def. prop.} \\
 & 2 & &
 \vdash \neg\neg ff(P) & &
 \text{Glivenko's thm.} \\
 & 3 & &
 ff(P) \rightarrow fff(P) \vdashl ff(P)  & &
 \text{cut on 1,2} \\
\end{align*}
The third claim is a consequence of the previous two:
\begin{align*}
 & \textbf{No.} & & \textbf{Claim} & & \textbf{Justification} \\
 & 1 & &
 \neg\neg f(P), f(P) \rightarrow ff(P) \vdashl f(P) & &
 \text{def. prop.} \\
 & 2 & &
 \vdash \neg\neg f(P) & &
 \text{Glivenko's thm.} \\
 & 3 & &
 f(P) \rightarrow ff(P) \vdashl f(P) & &
 \text{cut on 1,2} \\
 & 4 & &
 f(P), ff(P) \vdashl fff(P) & &
 \text{claim 1} \\
 & 5 & &
 f(P) \vdashl ff(P) \rightarrow fff(P) & &
 \text{$\rightarrow R$ on 4} \\
 & 6 & &
 ff(P) \rightarrow fff(P) \vdashl ff(P) & &
 \text{claim 2} \\
 & 7 & &
 f(P) \vdashl ff(P) & &
 \text{cut on 5,6} \\
 & 8 & &
 \vdashl f(P) \rightarrow ff(P) & &
 \text{$\rightarrow R$ on 7} \\
 & 9 & &
 \vdashl f(P) & &
 \text{cut on 3,8} \\
\end{align*}
Note that we can restate this as $\vdashl t(P) \vee \neg t(P)$ by expanding definitions.
\end{proof}
\end{lemma}

\begin{theorem}\label{thm:polacik:mainresult}
We can infer $\neg\neg P \vdashl P$ for any formula $P$.
\begin{proof}
Immediate from Lemma~\ref{lemma:polacik:three} and the defining property of $t$.
\end{proof}
\end{theorem}

\begin{body}
As in the previous cases, Theorem~\ref{thm:polacik:undefinable} follows immediately from Theorem~\ref{thm:polacik:mainresult}: if intuitionistic propositional logic could define the connective $\bullet(P)$, then by Theorem~\ref{thm:general:auxiliaryequiv} its auxiliary term would satisfy the defining properties above, and therefore one could construct a proof of $\neg\neg P \vdash P$ for any $P$. Consequently, one cannot define Po\l acik's connective by a quantifier-free schema.
\end{body}

\subsection{Consequences}

\begin{body}
Analysis of the structural proof allows us on one hand to strengthen the result and obtain new non-definability results about connectives that share some structure with Po\l acik's connective, and on the other hand to obtain information about the flavors of logic which \textit{do} define the $\bullet(P)$ connective.
\end{body}

\begin{corollary}\label{cor:polacik:wlem}
The regular connective $C_\Phi(X)$ given by the formula $$(X \rightarrow (\neg Y \vee \neg\neg Y)) \rightarrow X$$ is not definable by a quantifier-free schema in intuitionistic propositional logic.
\begin{proof}
A \texttt{propquant} computation reveals that $\pite Y. (X \rightarrow (\neg Y \vee \neg\neg Y)) \rightarrow X$, just like the corresponding Pitts interpolant of Po\l acik's connective, coincides with $\neg \neg X$. From there, one can see that the argument of Theorem~\ref{thm:polacik:mainresult} applies without modification.
\end{proof}
\end{corollary}

\begin{corollary}\label{cor:polacik:disjunction}
Any intermediate logic $\mathcal{L}$ that provides a quantifier-free definition of the connective $\bullet(P)$ must prove
$$ \vdashl (\neg\neg X) \vee (\neg\neg X \rightarrow X). $$
Consequently, no such logic has the disjunction property.
\begin{proof}
Since we already calculated $\pite Y. (X \rightarrow (Y \vee \neg Y)) \rightarrow X$ as $\neg\neg X$ in Proposition~\ref{prop:polacik:interpolant}, we know that such a logic satisfies one of
\begin{enumerate}
    \item $\vdashl \neg\neg X \leftrightarrow \bullet(X)$
    \item $\vdashl X \leftrightarrow \bullet(X)$, or
    \item $\vdashl \bot \leftrightarrow \bullet(X)$.
\end{enumerate}
As before, the first case reduces to an application of Theorem~\ref{thm:polacik:mainresult}. The third case is trivial, which leaves only the case of $\vdashl X \leftrightarrow \bullet(X)$. To settle this final case, we now establish $\vdashl (\neg\neg X) \vee (\neg\neg X \rightarrow X)$ from the assumption $(X \rightarrow (Y \vee \neg Y)) \rightarrow X \vdashl X$ using an elementary, but surprisingly tricky argument. The key step involves making the substitutions
\begin{enumerate}
    \item $P \vee (P \rightarrow (Q \vee \neg Q))$ for $X$, and
    \item $Q$ for $Y$,
\end{enumerate}
in $X \rightarrow (Y \vee \neg Y) \vdashl X$ to obtain the monstrous sequent
$$ (P \vee (P \rightarrow (Q \vee \neg Q))) \rightarrow (Q \vee \neg Q)) \rightarrow (P \vee (P \rightarrow (Q \vee \neg Q))) \vdashl P \vee (P \rightarrow (Q \vee \neg Q)),$$
and carefully checking that the left hand side of the sequent above holds already as a tautology of intuitionistic propositional logic. This allows one to conclude the sequent $\vdashl P \vee (P \rightarrow (Q \vee \neg Q))$. From there, one quickly gets $\vdashl P \vee (P \rightarrow \neg\neg Q \rightarrow Q)$ by ordinary intuitionistic reasoning, and concludes by substituting first $\neg\neg X$ for $P$ and then $X$ for $Q$.
\end{proof}
\end{corollary}

\begin{body}
As before, the type-theoretic analogue of Theorem~\ref{thm:polacik:mainresult} follows immediately. See the file \texttt{Po\l acik.agda} of the repository for the full verification. We state the main result as Proposition~\ref{prop:polacik:univalent} below.
\end{body}

\begin{proposition}[in Agda]\label{prop:polacik:univalent}
Assume we have some $t: \mathtt{Prop} \rightarrow \mathtt{Prop}$ and an inhabitant of the following type:
$$ \Pi P: \mathtt{Prop}. ((P \Rightarrow \bot) \Rightarrow \bot) \Rightarrow (P \Rightarrow (t\:P \vee (t\:P \Rightarrow \bot))) \Rightarrow P.$$
Then we can construct an inhabitant of the type $\Pi P:\mathtt{Prop}. ((P \Rightarrow \bot) \Rightarrow \bot) \Rightarrow P$.
\end{proposition}

%%%%%%%%%%%%%%%%%%%%%%%%%%%%%%%%%%%%%%%%%%%%%%%%%%%%%%%%%%%%%%%%%%%%%%%%

\section{Future work}

\begin{body}
We have seen that the proof-theoretic method presented here can be used to settle quantifier-free nondefinability questions of interest (realizability disjunction) for regular connectives, and to improve known nondefinability results about unary connectives (Kreisel, Po\l acik). While we did not state the results of Section~\ref{sec:general} in full generality, it's clear that the method extends beyond the case of regular connectives with $\vee$-free interpolants. The relevant results are best developed on a case-by-case basis, as needed in specific applications. One could extend the applicability of the ideas even further, by replacing intuitionistic logic itself with certain super-intuitionistic logics, or with certain subsystems of intuitionistic logic. However, one must remember that the identity of the Pitts interpolants is closely tied to the logic under consideration: a formula that satisfies the definition of uniform interpolant in intuitionistic logic will generally not satisfy the same condition in KC, or even in the negation-free fragment of intuitionistic logic. Fortunately, in subsystems of intuitionistic logic, the uniform interpolants can frequently be computed from the Pitts interpolants themselves. For example, in the aforementioned negation-free fragment, the De Jongh-Zhao~\cite{dejongh-positive} positive part operator allows us to calculate the uniform interpolant from the result returned by \texttt{propquant} in linear time. Moreover, Iemhoff~\cite{iemhoff-uniform} recently related interpolation in both sub-structural and super-intuitionistic logics to the existence of so-called ``centered'' proof calculi for such logics. In principle, one could generalize the Coq formalization from which \texttt{propquant} was derived to work parametrically in an arbitrary such proof calculus, which would further expand the applicability of the method to a wide variety of sub-structural and modal logics, as well as to those super-intuitionistic logics which admit uniform interpolation.
\end{body}

\begin{body}
The fact that one can always choose the auxiliary formulae of Theorem~\ref{thm:general:auxiliaryequiv} in a quantifier-free way enables proof techniques that exploit the Rieger-Nishimura lattice in the one-variable case, and suggests that decidability results for the definability question may be attainable for certain classes of second-order connectives, in particular regular and $\forall$-regular connectives in one free variable. We leave these investigations for future work.
\end{body}

\section*{Acknowledgments}

\noindent We thank the anonymous referees for their careful reading of the manuscript, the suggestions for streamlining the presentation, and the informal summary of Proposition~\ref{prop:general:definability}. We thank Iris van der Giessen, Vincent Jackson and Christine Rizkallah for fruitful discussions, and Anupam Das for the suggestion to obtain the quantifier-free clause of Theorem~\ref{thm:general:auxiliaryequiv} directly from Pitts' translation.

\newpage

\bibliography{biblio}
\bibliographystyle{plainurl}

\end{document}